\documentclass[10pt]{amsart}
\usepackage[latin1]{inputenc}
\usepackage[T1]{fontenc}
\usepackage{indentfirst}
\usepackage{textcomp}
\usepackage[english]{babel}
\usepackage{amsmath,amsfonts,amssymb,amsopn,amscd,amsthm}
\usepackage{mathbbol}
\usepackage{mathrsfs}
\usepackage{graphicx}
\usepackage{pstricks,pst-plot}
\usepackage{xy}  \xyoption{all}
\usepackage{pb-diagram,pb-xy}
\usepackage{enumerate}

\begin{document}

\author{Pierre-Lo\"ic M\'eliot}
\title[Gaussian concentration of the $q$-characters of the {H}ecke algebras of type {A}]{Gaussian concentration of the $q$-characters of the {H}ecke algebras of type {A}}
\address{The Gaspard--Monge Institut of electronic and computer science,
University of Marne-La-Vall\'ee Paris-Est,
77454 Marne-la-Vall\'ee Cedex 2, France}
\email{meliot@phare.normalesup.org}

\begin{abstract}
We show that with respect to the $q$-Plancherel measure on partitions of size $n$, the irreducible characters of an Hecke algebra $\mathscr{H}_{q}(\mathfrak{S}_{n})$ are concentrated around the normalized trace of $\mathscr{H}_{q}(\mathfrak{S}_{n})$. More precisely, we prove that the deviations of the values of the $q$-characters $\chi^{\lambda}_{q}$ are asymptotically gaussian, and we give an explicit formula for the covariances of the limit normal laws. Our proof involves \'Sniady's theory of cumulants of observables of diagrams and a M\"obius inversion formula for additive class functions on symmetric groups.
\end{abstract}
\maketitle

\newcommand{\Z}{\mathbb{Z}}                
\newcommand{\N}{\mathbb{N}}                
\newcommand{\R}{\mathbb{R}}                     
\newcommand{\C}{\mathbb{C}}            
\newcommand{\IH}{\mathscr{H}}          
\newcommand{\sym}{\mathfrak{S}}      
\newcommand{\Comp}{\mathfrak{C}}   
\newcommand{\Part}{\mathfrak{P}}      
\newcommand{\proba}{\mathbb{P}}      
\newcommand{\esper}{\mathbb{E}}      
\newcommand{\GL}{\mathrm{GL}}       
\newcommand{\GB}{\mathrm{B}}         
\newcommand{\For}{\mathbb{F}}         
\newcommand{\card}{\mathrm{card}}   
\newcommand{\lle}{\left[\!\left[}              
\newcommand{\rre}{\right]\!\right]}    
\newcommand{\scal}[2]{\left\langle #1\vphantom{#2}\,\right |\left.#2 \vphantom{#1}\right\rangle}   
\newtheorem{theorem}{Theorem}
\newtheorem{proposition}[theorem]{Proposition}

If $G$ is a finite group and $V$ is a finite-dimensional complex linear representation of $G$, the decomposition in irreducible components $V=\bigoplus_{\lambda \in \widehat{G}} n_{\lambda}\,V_{\lambda}$ yields a probability measure on the set $\widehat{G}$ of isomorphism classes of irreducible representations of $G$:
$$\proba_{V}\!\left[\lambda \in \widehat{G}\right]=\frac{n_{\lambda}\,\dim V_{\lambda}}{\dim V}$$
The same goes for a complex linear representation of a semisimple finite-dimensional algebra $\mathscr{A}$: it provides a probability measure on $\widehat{\mathscr{A}}$. The Plancherel measures of the symmetric groups are well-known examples of such measures, see \cite{IO02}; in this case, the irreducible representations are labelled by partitions. We shall denote by $\sym_{n}$ the symmetric group of order $n$, by $\Part_{n}$ the set of integer partitions $\lambda$ of size $n$, and by $\chi^{\lambda}(k)$ the value of the normalized irreducible character with label $\lambda$ on a cycle of length $k$. If $\lambda$ is picked randomly according to the Plancherel measure $M_{n}$ associated to the left regular representation $\sym_{n} \curvearrowright \C\sym_{n}$, then 
$$\forall k \geq 2,\,\,\, n^{\frac{k}{2}}\,\chi^{\lambda}(k) \longrightarrow X_{k},$$
where the $X_{k}$'s are independant centered normal laws of variance $k$, see Theorem 6.1 in \cite{IO02}.\bigskip\bigskip

In this paper, we show that the same phenomenon occurs for the values of the characters of the Hecke algebras $\IH_{q}(\sym_{n})$, assuming that the irreducible characters are picked randomly with respect to the so-called $q$-Plancherel measures $M_{n,q}$. \medskip

\begin{theorem}\label{global}
In the following, $q$ is a real parameter in $\R_{+}^{*}\setminus \{1\}$, and we denote by $\chi^{\lambda}_{q}(k)$ the value of the normalized irreducible character $\chi^{\lambda}_{q}$ of the Hecke algebra $\IH_{q}(\sym_{n})$ on a ``cycle'' $T_{1}T_{2}\cdots T_{k-1}$, where the $T_{i}$'s  are the usual generators of $\IH_{q}(\sym_{n})$. As $n$ goes to infinity, if $\lambda$ is picked randomly according to the $q$-Plancherel measure $M_{n,q}$, then 
$$\forall k \geq 2,\,\,\,\sqrt{n}\,\,\chi^{\lambda}_{q}(k) \longrightarrow X_{q,k},$$
where $X_{q,k}$ is a centered normal law, and the arrow means that we have convergence in law. The covariances of the $X_{q,k}$'s are:
$$\mathrm{cov}(X_{q,k},X_{q,l})=(q-q^{2})^{k+l-3}\,(1-q^{2})\,\frac{\{k-1\}_{q}\,\{l-1\}_{q}}{\{k+l-1\}_{q}\,\{k+l-2\}_{q}\,\{k+l-3\}_{q}}$$
where $\{n\}_{q}$ is the $q$-analog of $n$, that is to say, $\frac{q^{n}-1}{q-1}$.
\end{theorem}\bigskip

\noindent The probability measures $M_{n,q}$ on partitions are the adequate quantizations of the Plancherel measures $M_{n}$ in the setting of Iwahori-Hecke algebras; we shall recall this in paragraph \ref{iwahorihecke}. These measures have already been studied in \cite{FM10}, but with a geometric point of view; namely, we were interested in the limit shape of the partitions under these probability laws.  In paragraphs \ref{newobs} and \ref{expectation}, we recall the arguments that lead to the gaussian concentration of $q$-characters; we already knew in \cite{FM10} that $q$-characters were asymptotically gaussian, but we did not know how to compute their actual limit laws. This computation is the true novelty in this paper, and it is done in paragraphs \ref{gaussian} and \ref{strange} by combining:\vspace{1mm}
\begin{itemize}
\item a formula for the second higher term of a product of two classes $\varSigma_{\mu}$ and $\varSigma_{\rho}$ in the Ivanov-Kerov algebra of partial permutations (\emph{cf.} \cite{IK99});\vspace{1mm}
\item a formula of Ram (\emph{cf.} \cite{Ram91}) that relates the $q$-characters of the Hecke algebra $\IH_{q}(\sym_{n})$ and the ``usual'' characters of the symmetric group $\sym_{n}$;\vspace{1mm}
\item and a remarkable identity that can be seen as a M\"obius inversion formula for additive class functions on symmetric groups.\vspace{1mm}
\end{itemize}
Besides, we make constant use of the algebra of polynomial functions on diagrams (see \cite[\S1-4]{IO02}), and of \'Sniady's theory of cumulants of such observables (\emph{cf.} \cite{Sni06}); our result is another evidence that they are extremely versatile tools in this setting of asymptotic representation theory.
\bigskip

\section{Hecke algebras of type A and their $q$-Plancherel measures}\label{iwahorihecke}
Let $q$ be a prime power, and $k=\For_{q}$ be the finite field with $q$ elements. We denote by $G=\GL(n,\For_{q})$ the group of $n \times n$ invertible matrices over $\For_{q}$, and by $B=\GB(n,\For_{q})$ the Borel subgroup of $G$ that consists in upper triangular matrices. The cosets in $G/B$ form a variety that parametrizes the complete flags in $(\For_{q})^{n}$, and there is a natural action of $G$ on this flag variety
$$g\cdot(hB)=(gh)B,$$
whence a complex linear representation of $\GL(n,\For_{q})$ on $\C[G/B]$. The irreducible components of this representation are the so-called unipotent modules $U_{\lambda}(\For_{q})$, and they are labelled by partitions $\lambda \in \Part_{n}$. Following ideas of E. Strahov, we defined in \cite{FM10} the $q$-Plancherel measure $M_{n,q}$ as the probability measure on partitions of size $n$ associated to the $\GL(n,\For_{q})$-module $\C[G/B]$:
$$\C[G/B]=\bigoplus_{\lambda \in \Part_{n}}n_{\lambda}\,U_{\lambda}(\For_{q}) \quad\Rightarrow\quad M_{n,q}(\lambda)=\frac{n_{\lambda}\,\dim U_{\lambda}(\For_{q})}{\card\,\GL(n,\For_{q})/\GB(n,\For_{q})}$$
The denominator is $\{n!\}_{q}=\{n\}_{q}\times\{n-1\}_{q} \times \cdots\times \{1\}_{q}$; in the numerator, $\dim U_{\lambda}(q)=D_{\lambda}(q)$ is the generic degree of label $\lambda$ and is a polynomial in $q$ (see \cite[Chapter 8]{GP00}), and it can be shown that $n_{\lambda}$ equals the dimension $\dim \lambda$ of the irreducible representation of $\sym_{n}$ of type $\lambda$. Indeed, the commutant of the action of $\GL(n,\For_{q})$ on $\C[G/B]$ is known since \cite{Iwa64} to be the Hecke algebra of $\sym_{n}$, that is to say, the complex algebra $\IH_{q}(\sym_{n})$ with generators $T_{1},\ldots,T_{n-1}$ and set of relations:
$$\begin{cases}&\forall i,\,\,\,(T_{i}-q)(T_{i}+1)=0\\
&\forall i,\,\,\,T_{i}T_{i+1}T_{i}=T_{i+1}T_{i}T_{i+1}\\
&\forall i,j,\,\,\,|i-j|\geq 2 \Rightarrow T_{i}T_{j}=T_{j}T_{i}
\end{cases}$$
It can be shown that $\IH_{q}(\sym_{n})$ admits for basis the $T_{\sigma}=T_{i_{1}}T_{i_{2}}\cdots T_{i_{r}}$, where $\sigma$ runs over $\sym_{n}$ and $\sigma=s_{i_{1}}s_{i_{2}}\cdots s_{i_{r}}$ is any reduced expression of $\sigma$ as a product of elementary transpositions $s_{i}=(i,i+1)$. Moreover, for generic $q$, the Iwahori-Hecke algebra $\IH_{q}(\sym_{n})$ is semisimple and has the same representation theory as the group algebra $\C\sym_{n}=\IH_{1}(\sym_{n})$. Then, if $V_{\lambda,q}$ denotes the irreducible representation of $\IH_{q}(\sym_{n})$ labelled by a partition $\lambda \in \Part_{n}$, one has the following decomposition of $\C[G/B]$ in irreducible $(\GL(n,\For_{q}),\IH_{q}(\sym_{n}))$-bimodules:
$$_{\GL(n,\For_{q})\curvearrowright}\C[G/B]_{\curvearrowleft \IH_{q}(\sym_{n})}=\bigoplus_{\lambda \in \Part_{n}} {}_{\GL(n,\For_{q})\curvearrowright}(U_{\lambda}(\For_{q}))\otimes_{\C} (V_{\lambda,q})_{\curvearrowleft \IH_{q}(\sym_{n})}$$
Consequently, $n_{\lambda}$ is indeed equal to $\dim \lambda$, and one can see $M_{n,q}$ as the probability measure associated to a representation of the Hecke algebra $\IH_{q}(\sym_{n})$. \bigskip\bigskip

More precisely, let us denote by $\tau_{q}$ the restriction to $\IH_{q}(\sym_{n})=\mathrm{End}_{\GL(n,\For_{q})}(\C[G/B])$ of the normalized trace of matrices in $\mathrm{End}_{\C}(\C[G/B])$. This symmetric trace is defined on the basis $(T_{\sigma})_{\sigma \in \sym_{n}}$ of the Hecke algebra by
$$\tau_{q}(T_{\sigma})=\begin{cases}1 & \text{if }\sigma =1,\\
0&\text{otherwise}.\end{cases}$$
Because of the aforementioned decomposition of the $(\GL(n,\For_{q}),\IH_{q}(\sym_{n}))$-bimodule $\C[G/B]$, if $\chi^{\lambda}_{q}$ denotes the normalized character of the irreducible representation $V_{\lambda,q}$ of $\IH_{q}(\sym_{n})$, then:
$$\tau_{q}=\sum_{\lambda \in \Part_{n}}M_{n,q}(\lambda)\,\chi^{\lambda}_{q}$$
So, the $q$-Plancherel measure $M_{n,q}(\lambda)$ is the weight of the $q$-character $\chi^{\lambda}_{q}$ in the canonical trace of the algebra $\IH_{q}(\sym_{n})$. With this latter definition, it makes sense to take $q$ in $\R_{+}^{*}$, and not only in the set of prime powers; in particular, if one sets $q=1$, one recovers the usual Plancherel measure $M_{n}$ of the symmetric group $\sym_{n}$, because $\IH_{1}(\sym_{n})=\C\sym_{n}$ and $\tau_{1}$ is the normalized character of the left regular representaion of the symmetric group.
\bigskip\bigskip

That said, we are interested here with the following problem. Let $q$ be a positive real number that is not equal to $1$. We suppose that $\lambda$ is chosen randomly among partitions of size $n$ according to the $q$-Plancherel measure $M_{n,q}$, and we fix an element $a \in \IH_{q}(\sym_{n})$, for instance a basis element $T_{\sigma}$. We then ask for the distribution of the random variable $\chi^{\lambda}_{q}(T_{\sigma})$. Because of the interpretation of $M_{n,q}(\lambda)$ as the weight of $\chi^{\lambda}_{q}$ in $\tau_{q}$, the expectation of this random variable is easy to compute:
$$\esper[\chi^{\lambda}_{q}(T_{\sigma})]=\tau_{q}(T_{\sigma})=\mathbb{1}_{\sigma=1}$$
We shall see that we have in fact convergence in probability of the $q$-characters towards their means, and with a gaussian concentration.\bigskip\bigskip

\section{Two bases of the algebra of observables of diagrams}\label{newobs}
A powerful tool in asymptotic representation theory of the symmetric groups is the algebra of observables of diagrams, also known as Kerov's algebra of polynomial functions on Young diagrams (\cite{KO94}). In the following, we present a graded basis of this algebra (the $\varSigma_{\mu}$'s), and a quantization of this basis (the $\varSigma_{\mu,q}$'s). We start with the symbols $\varSigma_{\mu}$ --- the so-called central characters --- and for now we consider them as conjugacy classes in the Ivanov-Kerov algebra of partial permutations $\mathscr{B}_{\infty}$, see \cite{IK99} for a precise definition, and also \cite[\S2 and \S4.3]{Sni06}. Hence, 
$$\varSigma_{\mu}=\sum_{a_{11}\neq a_{12} \neq \cdots \neq a_{r\mu_{r}} } (a_{11},a_{12},\ldots,a_{1\mu_{1}})(a_{21},\ldots,a_{2\mu_{2}})\cdots (a_{r1},\ldots,a_{r\mu_{r}}),$$
where the $a$'s run over the set of positive integers $\N^{*}$  (\emph{cf.} \cite[\S2.1]{Sni06}). If $n\geq 0$, then one can project such a symbol in the center of the symmetric group algebra $\C\sym_{n}$ by requiring that the $a$'s remain in the interval $\lle 1,n\rre$. One obtains a multiple of the conjugacy class of type $\mu1^{n-|\mu|}$ if $n \geq |\mu|$, and $0$ otherwise.\bigskip\bigskip

The projections $\mathrm{pr}_{n}:\mathscr{B}_{\infty} \to \C\sym_{n}$ described above form a separating family of morphisms of algebras, and consequently, the symbols $\varSigma_{\mu}$ generate inside $\mathscr{B}_{\infty}$ a commutative subalgebra $\mathscr{A}_{\infty}$. More precisely, if 
$$\varSigma_{\mu}=\!\!\!\!\!\!\sum_{a_{11}\neq \cdots \neq a_{r\mu_{r}} }\!\!\!\! (a_{11},\ldots,a_{1\mu_{1}}) \cdots (a_{r1},\ldots,a_{r\mu_{r}}) \quad\text{and}\quad \varSigma_{\nu}=\!\!\!\!\!\!\sum_{b_{11}\neq \cdots \neq b_{s\rho_{s}} }\!\!\!\! (b_{11},\ldots,b_{1\rho_{1}}) \cdots (b_{s1},\ldots,b_{s\rho_{s}})$$
then $\varSigma_{\mu}\varSigma_{\nu}=\sum_{M} \varSigma_{\rho(M)}$, where $M$ runs over the partial matchings of the set of indices of the $a$'s with the set of indices of the $b$'s, and $\rho(M)$ depends only on the partial matching $M$, see \cite[\S3.3]{FM10}. For instance, suppose that $\mu=(3)$ and $\rho=(2)$, and let us compute the product $\varSigma_{3}\varSigma_{2}$; the sets of indices are $\{1,2,3\}$ for the $a$'s and $\{1',2'\}$ for the $b$'s. \vspace{2mm}
\begin{enumerate}
\setcounter{enumi}{-1}
\item The empty partial matching corresponds to products of a $3$-cycle with a disjoint $2$-cycle, whence a contribution $\varSigma_{3,2}$.\vspace{2mm}
\item The six partial matchings of size $1$ correspond to products of a $3$-cycle $(a,b,c)$ with a $2$-cycle $(c,d)$, that is to say, to $4$-cycles $(a,b,c,d)$. Hence, these matchings give $6\varSigma_{4}$.\vspace{2mm}
\item The six partial matchings of size $2$ correspond to products of a $3$-cycle $(a,b,c)$ with a $2$-cycle $(b,c)$, that is to say, to disjoint products $(a,b)(c)$. So, these matchings contribute to $6\varSigma_{2,1}$.\vspace{2mm}
\end{enumerate}
Consequently, $\varSigma_{3}\varSigma_{2}=\varSigma_{3,2}+6\varSigma_{4}+6\varSigma_{2,1}$. In the general case, the size of a partition $\rho(M)$ that corresponds to a partial matching $M$ of size $|M|$ is $|\rho(M)|=|\mu|+|\nu|-|M|$, because this is the number of distinct indices once the matching is done. So, for all partitions $\mu$, $\nu$, the product $\varSigma_{\mu}\varSigma_{\nu}$ may be written as
$$\varSigma_{\mu\sqcup\nu}+\sum_{|M|=1}\varSigma_{\rho(M)}+\sum_{|\rho|\leq |\mu|+|\nu|-2}*\,\varSigma_{\rho},$$
where $\mu \sqcup \nu$ is the partition whose parts are those of $\mu$ and those of $\nu$. Let us precise the second term of this expansion. If $M$ is a partial matching of size $1$, it means that one identifies exactly one of the $a$'s in a cycle $(a_{i1},\ldots,a_{i\mu_{i}})$ with of the $b$'s in a cycle $(b_{j1},\ldots,b_{j\nu_{j}})$. The product of these cycles is therefore a $(\mu_{i}+\nu_{j}-1)$-cycle, and the remaining cycles stay disjoint and keep their respective sizes. If two parts $c=\mu_{i}$ and $d=\nu_{j}$ are fixed, there are $cd$ corresponding partial matchings of size $1$; as a consequence,
$$\varSigma_{\mu}\varSigma_{\nu}-\varSigma_{\mu\sqcup\nu}=\sum_{\substack{c \in \mu\\d\in \nu}} cd\, \varSigma_{(\mu\setminus\{c\})\sqcup(\nu \setminus \{d\})\sqcup (c+d-1)}+\sum_{|\rho|\leq |\mu|+|\nu|-2}*\,\varSigma_{\rho}.$$
This expansion will play a prominent role in the computation of the cumulants of the observables $\varSigma_{\mu}$ under $q$-Plancherel measure, see \S\ref{gaussian}.\bigskip\bigskip

From now on, we consider the elements of $\mathscr{A}_{\infty}$ as observables of diagrams, meaning that we evaluate them on partitions by following the rule 
$$\varSigma_{\mu}(\lambda)=\chi^{\lambda}(\mathrm{pr}_{n}(\varSigma_{\mu}))=\begin{cases}n^{\downarrow k}\,\chi^{\lambda}(\mu1^{n-k})&\text{if }|\mu|=k\leq n=|\lambda|,\\
0&\text{otherwise.}
\end{cases}$$ 
Here, $n^{\downarrow k}$ is the $k$-th falling factorial of $n$, and by $\chi^{\lambda}(\mu1^{n-k})$, we mean the value of the irreducible character $\chi^{\lambda}$ on a permutation $\sigma$ of cycle type $\mu1^{n-k}$. Notice that these evaluations of elements in $\mathscr{A}_{\infty}$ are compatible with the product of partial permutations: 
$$\forall \lambda,\,\,\forall a,b \in \mathscr{A}_{\infty},\,\,[ab](\lambda)=a(\lambda)\,b(\lambda)$$
Indeed, if $\chi^{\lambda}$ is a normalized irreducible character on $\C\sym_{n}$ and if $a$ and $b$ are two elements in the center of the symmetric group algebra, then $\chi^{\lambda}(ab)=\chi^{\lambda}(a)\,\chi^{\lambda}(b)$, because $\chi^{\lambda}(x)$ is the eigenvalue of the action of a central element $x$ on the irreducible representation space $V_{\lambda}$. That said, one 
introduces $q$-symbols $\varSigma_{\mu,q}$ that are defined by the evaluations:
$$\varSigma_{\mu,q}(\lambda)=\begin{cases}n^{\downarrow k}\,\chi^{\lambda}_{q}(\mu1^{n-k})&\text{if }|\mu|=k\leq n=|\lambda|,\\
0&\text{otherwise.}
\end{cases}$$ 
Here, $\chi^{\lambda}_{q}(\mu1^{n-k})$ is the value of the irreducible character $\chi^{\lambda}_{q}$ of the Hecke algebra $\IH_{q}(\sym_{n})$ on a basis element $T_{\sigma}$, with $\sigma$ of minimal length among permutations with cycle type $\mu1^{n-k}$. A formula due to A. Ram (\emph{cf.} \cite[Theorem 5.4]{Ram91} and \cite[Proposition 10]{FM10}) relates the $q$-symbols $\varSigma_{\mu,q}$ to the regular symbols $\varSigma_{\mu}$:
\begin{align*}
\forall \rho \in \Part_{k},\,\,\,(q-1)^{\ell(\rho)}\,\varSigma_{\rho,q}(\lambda)&=\sum_{\nu \in \Part_k} \frac{q^\nu-1}{z_\nu}\scal{p_\nu}{h_\rho}\,\varSigma_\nu(\lambda) \\
\forall \rho \in \Part_{k},\,\,\,(q^\rho-1)\,\varSigma_\rho(\lambda)&=\sum_{\nu \in \Part_k} (q-1)^{\ell(\nu)}\scal{m_\nu}{p_\rho}\,\varSigma_{\nu,q}(\lambda)\end{align*}
with\footnote{In the following, we use the same notation for $1-q^{\nu}=\prod_{i=1}^{\ell(\nu)}(1-q^{\nu_{i}})$. In particular, $1-q^{a,b}=(1-q^{a})\,(1-q^{b})$. } $q^{\nu}-1=\prod_{i=1}^{\ell(\nu)}(q^{\nu_{i}}-1)$ and $z_{\nu}=\scal{p_{\nu}}{p_{\nu}}=k! /\card\, C_{\nu}$. Here, $\scal{\cdot}{\cdot}$ is the usual scalar product in the algebra $\Lambda$ of symmetric functions, the $p$'s are the power sums, the $m$'s are the monomial symmetric functions and the $h$'s are the complete symmetric functions (see \cite{Mac95}). These formulae allow to treat the $\varSigma_{\rho,q}$ as elements of $\mathscr{A}_{\infty}$, because they are linear combinations of observables $\varSigma_{\nu}$. Thus, one has a commutative algebra of observables of diagrams with graded basis $(\varSigma_{\mu})_{\mu}$, and the deformation of the group algebras $\C\sym_{n}$ into the Hecke algebras $\IH_{q}(\sym_{n})$ yields a quantization $(\varSigma_{\mu,q})_{\mu}$ of the graded basis.\bigskip\bigskip

\section{Expectations and cumulants of observables}\label{expectation}
If $a$ is an observable of diagrams, we denote by $\esper[a]=M_{n,q}[a]$ the expectation of the random variable $a(\lambda)$ under the $q$-Plancherel measure $M_{n,q}$. The expectations of the $q$-symbols are easy to compute:
$$\esper[\varSigma_{\mu,q}]=n^{\downarrow |\mu|}\,\esper[\chi^{\lambda}_{q}(\mu1^{n-|\mu|})]=n^{\downarrow|\mu|} \,\mathbb{1}_{\mu=1^{k}} $$
with $n^{\downarrow |\mu|}=0$ if $|\mu|>n$. Now, by using the change of basis formula between characters and $q$-characters, one can also compute easily the expectations of the symbols $\varSigma_{\mu}$:
$$\esper[\varSigma_{\mu}]=\frac{1}{(q^{\mu}-1)}\sum_{\nu \in \Part_{k}}(q-1)^{\ell(\nu)}\scal{m_{\nu}}{p_{\mu}}\,\esper[\varSigma_{\nu,q}]=\frac{(q-1)^{|\mu|}}{q^{\mu}-1}\scal{m_{1^{|\mu|}}}{p_{\mu}}n^{\downarrow |\mu|}=\frac{(1-q)^{|\mu|}}{1-q^{\mu}}\,n^{\downarrow |\mu|}$$
In particular, $\esper[\varSigma_{\rho}]$ is always a $O(n^{|\rho|})$. From this, we deduce the following convergences in probability:
$$\forall \mu,\,\,\,\frac{\varSigma_{\mu}(\lambda)}{n^{|\mu|}} \longrightarrow_{M_{n,q}} \frac{(1-q)^{|\mu|}}{1-q^{\mu}}$$
Indeed, the expectation of the left-hand side is asymptotically equal to the right-hand side, and the variance of the left-hand side is
\begin{align*}\frac{\big(\esper[\varSigma_{\mu}^{2}]-\esper[\varSigma_{\mu}]^{2}\big)}{n^{2|\mu|}}&=\frac{\big(\esper[\varSigma_{\mu\sqcup\mu}]-\esper[\varSigma_{\mu}]^{2}\big)}{n^{2|\mu|}}+\frac{\esper[\text{some linear combination of }\varSigma_{\rho}\text{'s with }|\rho|\leq 2|\mu|-1]}{n^{2|\mu|}}\\
&=\frac{(1-q)^{2|\mu|}}{1-q^{\mu\sqcup \mu}}\,\frac{(n^{\downarrow 2|\mu|}-n^{\downarrow |\mu|}\,n^{\downarrow |\mu|})}{n^{2|\mu|}}+O(n^{-1})=O(n^{-1}),\end{align*}
whence the result by Bienaym\'e-Chebyshev inequality. Then, if one uses again the change of basis formula, one sees that $\varSigma_{\mu,q}(\lambda)/n^{|\mu|}$ converges in probability towards its asymptotic mean value, that is to say, $\mathbb{1}_{\mu=1^{k}}$. Since $\varSigma_{\mu,q}(\lambda)/n^{|\mu|}$ is essentially $\chi^{\lambda}_{q}(\mu1^{n-|\mu|})$, one has therefore proved: \medskip

\begin{proposition}
Let us denote by $T_{\mu}$ the basis element of Hecke algebras associated to the permutation 
$$\sigma_{\mu}=(1,2\ldots,\mu_{1})(\mu_{1}+1,\ldots,\mu_{1}+\mu_{2})\cdots (\mu_{1}+\cdots+\mu_{r-1}+1,\ldots,|\mu|)$$
that is of minimal length in its conjugacy class. When $n$ goes to infinity, $\chi^{\lambda}_{q}(T_{\mu})$ converges in probability towards the trace $\tau(T_{\mu})$. In particular, for all $k \geq 2$, $\chi^{\lambda}_{q}(k) \to 0$.
\end{proposition}
\bigskip\bigskip

Now, to prove the gaussian deviation of (rescaled) observables of diagrams $X_{i}$, one will consider their joint cumulants $k(X_{i_{1}},\ldots,X_{i_{r}})$ that are defined recursively by the equations:
$$\esper[X_{1}\cdots X_{r}]=\sum_{\pi \text{ set partition of } \lle 1,r \rre} k(X_{i \in \pi_{1}})\,k(X_{i \in \pi_{2}})\,\cdots \, k(X_{i \in \pi_{l}})$$
For instance, $k(X)=\esper[X]$ and $k(X,Y)$ is the covariance $\esper[XY]-\esper[X]\esper[Y]$. A gaussian vector $(X_{i})_{i \in I}$ is characterized by the fact that all cumulants $k(X_{i_{1}},\ldots,X_{i_{r}})$ of order $r$ higher than $3$ are equal to $0$. Then, the cumulants of order $2$ give the covariance matrix of the gaussian vector. In the setting of representation theory, the idea to use joint cumulants of observables of partitions in order to highlight phenomena of gaussian concentration is due to P. \'Sniady, \emph{cf.} \cite{Sni06}. For our purpose, we shall only need the following result:\medskip 

\begin{proposition}
If $a=\sum_{\rho} c_{\rho}\,\varSigma_{\rho}$ is an observable of diagrams in $\mathscr{A}_{\infty}$, the degree of $a$ is the highest size $|\rho|$ such that $c_{\rho}\neq 0$. Then, we have the following estimation for the joint cumulants of observables of diagrams under the $q$-Plancherel measure:
$$k(a_{1},\ldots,a_{r})=O(n^{\deg a_{1}+\deg a_{2}+\cdots+\deg a_{r}-r+1})$$ 
\end{proposition}\bigskip

\noindent This proposition is entirely proved in \cite[Lemma 16]{FM10}. We won't reproduce the proof here, but we can give an account of the main arguments:\vspace{1mm}
\begin{itemize}
\item One can define disjoint cumulants of observables $k^{\bullet}(a_{1},\ldots,a_{r})$ that correspond to another product in $\mathscr{A}_{\infty}$, namely, the ``disjoint product'' $\varSigma_{\mu}\bullet \varSigma_{\nu}=\varSigma_{\mu\sqcup \nu}$, see \cite[\S2.2]{Sni06}. For these disjoint cumulants, the same estimate
$$k^{\bullet}(a_{1},\ldots,a_{r})=O(n^{\deg a_{1}+\deg a_{2}+\cdots+\deg a_{r}-r+1})$$ 
holds, and this is easy to prove just by looking at the expectations of the symbols $\varSigma_{\mu}$.\vspace{2mm}
\item The (standard and disjoint) cumulants are multilinear and have a good behavior with respect to products of observables. On the other hand, disjoint cumulants and standard cumulants of observables are related one to another by a conditioning technique due to D. Brillinger, see \cite[\S4.4]{Sni06}. For these reasons, the problem is eventually reduced to the computation of the degree of the so-called identity cumulants $k^{\mathrm{id}}(\varSigma_{i_{1}},\ldots,\varSigma_{i_{r}})$. These observables of diagrams are defined recursively by the equations:
$$X_{1}\cdots X_{r}=\sum_{\pi \text{ set partition of } \lle 1,r \rre} k^{\mathrm{id}}(X_{i \in \pi_{1}})\bullet k^{\mathrm{id}}(X_{i \in \pi_{2}})\bullet\cdots \bullet k^{\mathrm{id}}(X_{i \in \pi_{l}})$$
and one has to show that $k^{\mathrm{id}}(\varSigma_{i_{1}},\ldots,\varSigma_{i_{r}})$ has degree less than $i_{1}+i_{2}+\cdots+i_{r}-r+1$.\vspace{2mm}
\item Finally, one can give a simple combinatorial interpretation of the observables $k^{\mathrm{id}}(\varSigma_{i_{1}},\ldots,\varSigma_{i_{r}})$. For any family $(a_{kl})_{1\leq l\leq i_{k},1\leq k \leq r}$ such that $a_{kl}\neq a_{kl'}$ for all $k \in \lle 1,r\rre$ and for all $l,l'$, let us consider the following relation on the integers $k \in \lle 1,r\rre$:
$$k_{1} \sim k_{2} \iff \exists l_{1}, l_{2},\,\,a_{k_{1}l_{1}}=a_{k_{2}l_{2}}$$
The relation $\sim$ can be completed in an equivalence relation on $\lle 1,r\rre$, and we will denote by $\pi((a_{kl})_{k,l})$ the associated set partition of $\lle 1,r\rre$. Then, it can be shown by recursion that:
$$k^{\mathrm{id}}(\varSigma_{i_{1}},\ldots,\varSigma_{i_{r}})=\sum_{(a_{kl})_{k,l}\,\,|\,\,\pi((a_{kl})_{k,l}=\{\lle 1,r\rre\}} (a_{11},\ldots,a_{1i_{1}})\cdots(a_{r1},\ldots,a_{ri_{r}})$$
In other words, in comparaison with $\varSigma_{i_{1}}\cdots\varSigma_{i_{r}}$, the sum is reduced to transitive products of cycles. Since this condition implies that there is at least $r-1$ identities between the $a$'s, there is at most $i_{1}+i_{2}+\cdots+i_{r}-r+1$ distinct $a$'s in the cycle decomposition of a partial permutation appearing in $k^{\mathrm{id}}(\varSigma_{i_{1}},\ldots,\varSigma_{i_{r}})$. Therefore, this identity cumulant is indeed a linear combination of $\varSigma_{\rho}$'s with $|\rho|\leq i_{1}+i_{2}+\cdots+i_{r}-r+1$.\vspace{2mm}
\end{itemize}
Again, we refer to \cite{FM10} for a complete and detailed proof of our claim; notice that it is not the same result as in \cite{Sni06}, because one uses a different gradation on the algebra $\mathscr{A}_{\infty}$.\bigskip\bigskip

\section{Gaussian deviation of characters and $q$-characters}\label{gaussian}
We now introduce the rescaled deviations of characters and $q$-characters:
$$Z_{\mu,n,q}= \sqrt{n}\left(\frac{\varSigma_{\mu}(\lambda)}{n^{|\mu|}}-\frac{(1-q)^{|\mu|}}{1-q^{\mu}}\right)\qquad;\qquad W_{k\geq 2,n,q}=\sqrt{n}\,\frac{\varSigma_{k,q}(\lambda)}{n^{k}}\simeq \sqrt{n}\,\chi^{\lambda}_{q}(k)$$
By construction, all these observables are centered random variables. Now, given partitions $\mu^{(1)},\ldots,\mu^{(r)}$ with $r\geq 3$, the joint cumulant of the $Z_{\mu^{(i)},n,q}$'s is:
\begin{align*}
k(Z_{\mu^{(1)},n,q},\ldots,Z_{\mu^{(r)},n,q})&=\frac{1}{n^{|\mu^{(1)}|-1/2}\cdots n^{|\mu^{(r)}|-1/2}}\,k(\varSigma_{\mu^{(1)}},\ldots,\varSigma_{\mu^{(r)}})\\
&=\frac{O(n^{|\mu_{1}|+\cdots+|\mu_{r}|-r+1})}{n^{|\mu_{1}|+\cdots+|\mu_{r}|-r/2}}=O\big(n^{1-r/2}\big)\to 0
\end{align*}
Consequently, the family of random variables $(Z_{\mu,n,q})_{\mu}$ converges in finite-dimensional laws towards a centered gaussian family $(Z_{\mu,\infty,q})_{\mu}$. Let us compute the covariances of the $Z_{\mu,\infty,q}$'s, that is to say, the limits of the cumulants
$k(Z_{\mu,n,q},Z_{\nu,n,q})$. We decompose such a cumulant in two parts:
\begin{align*}n^{|\mu|+|\nu|-1}\,k(Z_{\mu,n,q},Z_{\nu,n,q})&=\esper[\varSigma_{\mu}\varSigma_{\nu}]-\esper[\varSigma_{\mu}]\,\esper[\varSigma_{\nu}]\\
&=\big(\esper[\varSigma_{\mu}\varSigma_{\nu}]-\esper[\varSigma_{\mu\sqcup\nu}]\big)+\big(\esper[\varSigma_{\mu\sqcup\nu}]-\esper[\varSigma_{\mu}]\,\esper[\varSigma_{\nu}]\big)\end{align*}
The second part is simply $\frac{(1-q)^{|\mu|+|\nu|}}{1-q^{\mu\sqcup \nu}}\,(n^{\downarrow |\mu|+|\nu|}-n^{\downarrow |\mu|}\,n^{\downarrow |\nu|})$, and the asymptotic expansion of a falling factorial is 
$$n^{\downarrow k}=n^{k}-\frac{k(k-1)}{2}n^{k-1}+O(n^{k-2}).$$
Consequently, the second part amounts to $-|\mu|\,|\nu|\,\frac{(1-q)^{|\mu|+|\nu|}}{1-q^{\mu\sqcup \nu}}\,\,n^{|\mu|+|\nu|-1}+O(n^{|\mu|+|\nu|-2})$. As for the first part, we use the expansion of $\varSigma_{\mu}\varSigma_{\nu}-\varSigma_{\mu \sqcup \nu}$ presented in paragraph \ref{newobs}:
\begin{align*}
\esper[\varSigma_{\mu}\varSigma_{\nu}-\varSigma_{\mu \sqcup \nu}]&=\sum_{\substack{c \in \mu\\d\in \nu}} cd\, \esper[\varSigma_{(\mu\setminus\{c\})\sqcup(\nu \setminus \{d\})\sqcup (c+d-1)}]+\sum_{|\rho|\leq |\mu|+|\nu|-2}*\,\esper[\varSigma_{\rho}]\\
&=\frac{(1-q)^{|\mu|+|\nu|}}{1-q^{\mu\sqcup \nu}}\left(\sum_{\substack{c \in \mu\\d\in \nu}} cd\,\frac{(1-q^{c})(1-q^{d})}{(1-q)(1-q^{c+d-1})}\right)n^{|\mu|+|\nu|-1}+O(n^{|\mu|+|\nu|-2})
\end{align*}
Hence, by decomposing $-|\mu|\,|\nu|$ in $-\sum_{c,d}cd$, one obtains the following estimate:
\begin{align*}
k(Z_{\mu,n,q},Z_{\nu,n,q})&=\frac{(1-q)^{|\mu|+|\nu|}}{1-q^{\mu\sqcup \nu}}\sum_{\substack{c \in \mu\\d\in \nu}} cd\left(\frac{(1-q^{c})(1-q^{d})}{(1-q)(1-q^{c+d-1})}-1\right)+O(n^{-1})\\
&=q\,\frac{(1-q)^{|\mu|+|\nu|}}{1-q^{\mu\sqcup \nu}}\sum_{\substack{c \in \mu\\d\in \nu}} cd\left(\frac{(1-q^{c-1})(1-q^{d-1})}{(1-q)(1-q^{c+d-1})}\right)+O(n^{-1})
\end{align*}
We have then computed the covariances of the limit normal laws $Z_{\mu,\infty,q}$. Now, by using the change of basis formula between characters and $q$-characters, one sees that each variable $W_{k,n,q}$ is a linear combination of variables $Z_{\mu,n,q}$, and these latter variables converge jointly towards a centered gaussian vector. Hence, the $W_{k,n,q}$'s converge jointly towards a centered gaussian vector $(X_{k,q})_{k \geq 2}$, and the limit covariances are given by the following formula:
\begin{align*}\mathrm{cov}(X_{k,q},X_{l,q})&=\frac{1}{(q-1)^{2}}\sum_{\substack{\mu \in \Part_{k}\\ \nu \in \Part_{l}}} \frac{\scal{h_{k}}{p_{\mu}}\,\scal{h_{l}}{p_{\nu}}  }{z_{\mu}\,z_{\nu}}\,(q^{\mu \sqcup \nu}-1)\,k(Z_{\mu,\infty,q},Z_{\nu,\infty,q}) \\
&=q\,(1-q)^{k+l-3}\,\sum_{\substack{\mu \in \Part_{k}\\ \nu \in \Part_{l}}} \frac{(-1)^{\ell(\mu)+\ell(\nu)}}{z_{\mu}\,z_{\nu}}\left(\sum_{\substack{c \in \mu\\ d \in \nu}} cd\,\frac{1-q^{c-1,d-1}}{1-q^{c+d-1}}\right) \end{align*}
Indeed, the scalar products $\scal{h_{k}}{p_{\mu}}$  and $\scal{h_{l}}{p_{\nu}}$ are all equal to $1$ by the Frobenius-Schur formula, and the signs $(-1)^{\ell(\cdot)}$ come from the simplifications $(q^{(\cdot)}-1)/(1-q^{(\cdot)})$. To summarize, we have shown for now the following:\medskip

\begin{proposition}
When $n$ goes to infinity, the random sequence $(\sqrt{n}\,\chi^{\lambda}_{q}(k))_{k \geq 2}$ converges in finite-dimen\-sio\-nal laws towards a centered gaussian process $(W_{X,q})_{k \geq 2}$ whose covariances are given by the formula above.
\end{proposition}\bigskip

\noindent It remains to be seen that the really complex formula for covariances can in fact be reduced to a sum of only four rational fractions in $q$, and that this sum admits the simple (factorized) expression given in the statement of Theorem \ref{global}. This last part of our reasoning is quite a funny calculation.
\bigskip\bigskip

\section{A M\"obius inversion formula for additive class functions}\label{strange}
The reduction of the formula for covariances relies mainly on the following trick. Let $f$ be any function on the set of positives integers, and denote by $F$ the additive class function on $\sym_{n}$:
$$F(\sigma)=\sum_{c \text{ cycle of }\sigma}f(|c|)$$
For instance, if $n=7$ and $\sigma=(1,4,5)(2,3)(6,7)$, then $F(\sigma)=f(3)+2f(2)$. If $\lambda$ is a partition of size $n$, we also denote by $F(\lambda)$ the sum of the $f(\lambda_{i})$'s.
\begin{proposition}\label{stupidlemma}
For $n$ greater than or equal to $2$,
$$\frac{1}{n!}\sum_{\sigma\in \sym_{n}} (-1)^{\text{number of cycles of }\sigma}\,F(\sigma)=\frac{f(n-1)}{n-1}-\frac{f(n)}{n}.$$
\end{proposition}\bigskip

\noindent If one collects the permutations according to their cycle type, one sees that Proposition \ref{stupidlemma} implies the following: 
$$\sum_{\lambda \in \Part_{n}} \frac{(-1)^{\ell(\lambda)}}{z_{\lambda}}F(\lambda) = \frac{f(n-1)}{n-1}-\frac{f(n)}{n}$$
Proposition \ref{stupidlemma} is really a M\"obius inversion formula if we translate it for set partitions of $\lle 1,n\rre$, and gather the permutations according to their orbits. However, the simplest proof that one can give for Proposition \ref{stupidlemma} does not use this interpretation, and is by induction on $n$. The result is clearly true for $n=2$, and if it is true up to rank $n-1\geq 2$, then the sum over permutations $\sigma \in \sym_{n}$ may be decomposed in the following contributions:\vspace{1mm}
\begin{enumerate}
\item If $\sigma$ is one of the $(n-1)!$ cycles of length $n$, then $F(\sigma)=f(n)$, whence a contribution
$$A_{n}=-\frac{(n-1)!}{n!}\,f(n)=-\frac{f(n)}{n}.$$\vspace{1mm}
\item If $\sigma$ is one of the $(n-1)!$ cycles of length $n-1$ that move the integer $1$, then $F(\sigma)=f(1)+f(n-1)$, whence a contribution
$$A_{n-1}=\frac{(n-1)!}{n!}\,(f(1)+f(n-1))=\frac{f(1)}{n}+\frac{f(n-1)}{n}.$$\vspace{1mm}
\item Otherwise, the integer $1$ is in a cycle of $\sigma$ of length $k \in \lle 1,n-2\rre$. Once such a cycle has been chosen within the $(n-1)^{\downarrow k-1}$ possible $k$-cycles that move the integer $1$, the product $\tau$ of the remaining cycles of $\sigma$ may be considered as an element of $\sym_{n-k}$, so:
\begin{align*}A_{k}&=\frac{(n-1)^{\downarrow k-1}}{n!}\sum_{\tau \in \sym_{n-k}} (-1)^{1+\text{number of cycles of }\tau}\left(F(\tau)+f(k)\right)\\
&=\frac{1}{n}\left(-F(n-k)-\frac{f(k)}{(n-k)!}\,\sum_{\tau \in \sym_{n-k}} (-1)^{\text{number of cycles of }\tau}\right)\\
&=\frac{1}{n}\left(\frac{f(n-k)}{n-k}-\frac{f(n-k-1)}{n-k-1}\right)
\end{align*}
Indeed, the induction hypothesis holds for $n-k \in \lle 2,n-1\rre$, and since $n-k \geq 2$, the sum of signs is equal to zero, because even and odd permutations of size $N\geq 2$ have the same cardinality.\vspace{1mm}
\end{enumerate}
By summing all contributions, one obtains $\frac{f(n-1)}{n-1}-\frac{f(n)}{n}$, so the result is true up to rank $n$; we have then proved Proposition \ref{stupidlemma}.
\bigskip\bigskip

Let us come back to the computation of covariances $\mathrm{cov}(X_{k,q},X_{l,q})$. We fix a partition $\mu \in \Part_{k}$, and a part $c$ of $\mu$. If $f(n)=cn\,\frac{1-q^{c-1,n-1}}{1-q^{c+n-1}}$, then one has by Proposition \ref{stupidlemma} the following simplification:
$$\sum_{\nu\in \Part_{l}}\frac{(-1)^{\ell(\nu)}}{z_{\nu}}\sum_{d \in \nu}cd\,\frac{1-q^{c-1,d-1}}{1-q^{c+d-1}}=\frac{f(l-1)}{l-1}-\frac{f(l)}{l}=c\left(\frac{1-q^{c-1,l-2}}{1-q^{c+l-2}}-\frac{1-q^{c-1,l-1}}{1-q^{c+l-1}}\right)$$
Consequently, if $f_{1}(m)=m\,\frac{1-q^{m-1,l-1}}{1-q^{m+l-1}}$ and $f_{2}(m)=m\,\frac{1-q^{m-1,l-2}}{1-q^{m+l-2}}$, and if $F_{1}$ and $F_{2}$ are the corresponding additive class functions as in the statement of Proposition \ref{stupidlemma}, then:
$$\mathrm{cov}(X_{k,q},X_{l,q})=q\,(1-q)^{k+l-3}\,\left(\sum_{\mu \in \Part_{k}}\frac{(-1)^{\ell(\mu)}}{z_{\mu}}\,F_{2}(\mu)-\sum_{\mu \in \Part_{k}}\frac{(-1)^{\ell(\mu)}}{z_{\mu}}\,F_{1}(\mu)\right)$$
We use again Proposition \ref{stupidlemma} and we reduce the previous expression to:
\begin{align*}
\mathrm{cov}(X_{k,q},X_{l,q})&=q\,(1-q)^{k+l-3}\,\left(\frac{f_{2}(k-1)}{k-1}-\frac{f_{2}(k)}{k}-\frac{f_{1}(k-1)}{k-1}+\frac{f_{1}(k)}{k}\right)\\
&=q\,(1-q)^{k+l-3}\,\left(\frac{1-q^{k-2,l-2}}{1-q^{k+l-3}}-\frac{1-q^{k-1,l-2}}{1-q^{k+l-2}}-\frac{1-q^{k-2,l-1}}{1-q^{k+l-2}}+\frac{1-q^{k-1,l-1}}{1-q^{k+l-1}}\right)\\
&=\frac{q^{k-1}\,(1-q)^{k+l-2}\,(1-q^{l-1})}{1-q^{k+l-2}}\,\left(\frac{1-q^{l}}{1-q^{k+l-1}}-\frac{1-q^{l-2}}{1-q^{k+l-3}}\right)\\
&=\frac{q^{k+l-3}\,(1-q)^{k+l-2}\,(1-q^{l-1})\,(1-q^{k-1})\,(1-q^{2})}{(1-q^{k+l-1})\,(1-q^{k+l-2})\,(1-q^{k+l-3})}\\
&=(q-q^{2})^{k+l-3}\,(1-q^{2})\,\frac{\{k-1\}_{q}\,\{l-1\}_{q}}{\{k+l-1\}_{q}\,\{k+l-2\}_{q}\,\{k+l-3\}_{q}}
\end{align*}
This ends the proof of Theorem \ref{global}. Notice that when $q<1$, two $q$-characters $\chi^{\lambda}_{q}(k)$ and $\chi^{\lambda}_{q}(l)$ are always asymptotically positively correlated, whereas when $q>1$, it depends on the parity of $k+l$.\bigskip\bigskip

\bibliographystyle{alpha}
\bibliography{qcharacter}

\begin{thebibliography}{Ram91}

\bibitem[FM10]{FM10}
V.~F\'eray and P.-L. M\'eliot.
\newblock Asymptotics of q-{P}lancherel measures.
\newblock \texttt{arXiv:1001.2180v1 [math.RT]}, 2010.

\bibitem[GP00]{GP00}
M.~Geck and G.~Pfeiffer.
\newblock {\em Characters of Finite Coxeter Groups and Iwahori-Hecke Algebras},
  volume~21 of {\em London Mathematical Society Monographs}.
\newblock Oxford University Press, 2000.

\bibitem[IK99]{IK99}
V.~Ivanov and S.~Kerov.
\newblock The algebra of conjugacy classes in symmetric groups, and partial
  permutations.
\newblock In {\em Representation Theory, Dynamical Systems, Combinatorial and
  Algorithmical Methods III}, volume 256 of {\em Zapiski Nauchnyh Seminarov
  POMI}, pages 95--120, 1999.

\bibitem[IO02]{IO02}
V.~Ivanov and G.~Olshanski.
\newblock Kerov's central limit theorem for the {P}lancherel measure on {Y}oung
  diagrams.
\newblock In {\em Symmetric Functions 2001: Surveys of Developments and
  Perspectives}, volume~74 of {\em NATO Science Series II. Mathematics, Physics
  and Chemistry}, pages 93--151, 2002.

\bibitem[Iwa64]{Iwa64}
N.~Iwahori.
\newblock On the structure of the {H}ecke ring of a {C}hevalley group over a
  finite field.
\newblock {\em J. Faculty Science Tokyo University}, 10:215--236, 1964.

\bibitem[KO94]{KO94}
S.~Kerov and G.~Olshanski.
\newblock Polynomial functions on the set of {Y}oung diagrams.
\newblock {\em Comptes Rend. Acad. Sci. Paris, S\'erie I}, 319:121--126, 1994.

\bibitem[Mac95]{Mac95}
I.~G. Macdonald.
\newblock {\em Symmetric functions and {H}all polynomials}.
\newblock Oxford Mathematical Monographs. Oxford University Press, 2nd edition,
  1995.

\bibitem[Ram91]{Ram91}
A.~Ram.
\newblock A {F}robenius formula for the characters of the {H}ecke algebras.
\newblock {\em Invent. Math.}, 106:461--488, 1991.

\bibitem[\'S06]{Sni06}
P.~\'Sniady.
\newblock Gaussian fluctuations of characters of symmetric groups and of
  {Y}oung diagrams.
\newblock {\em Probab. Theory and Related Fields}, 136(2):263--297, 2006.

\end{thebibliography}

\end{document}